\documentclass[a4paper,12pt,reqno]{amsart}

\usepackage{a4wide}
\usepackage{amsmath}
\usepackage{amssymb}
\usepackage{theorem}
\usepackage{enumerate}
\usepackage{graphicx}
\usepackage{tikz}
\usetikzlibrary{arrows}
\newcommand{\gauss}[2]{\genfrac{[}{]}{0pt}{}{#1}{#2}}

\title[Lehmer's tridiagonal matrix]{The LU-decomposition of Lehmer's tridiagonal matrix}

\author[H.~Prodinger]{Helmut Prodinger}

\address{Helmut Prodinger,
	Mathematics Department, Stellenbosch University,
	7602 Stellenbosch, South Africa.}
\email{hproding@sun.ac.za}

\dedicatory{Dedication to George Andrews on the occasion of his 80\textsuperscript{th} birthday.}

\begin{document}

	\maketitle

\allowdisplaybreaks

\section{Introduction}

Ekhad and Zeilberger \cite{EZ} have unearthed Lehmer's \cite{Lehmer} tridiagonal $n\times n$ matrix $M=M(n)$ with entries
\begin{equation*}
M_{i,j}=\begin{cases}
1 & \text{if } i=j,\\
z^{1/2}q^{(i-1)/2}& \text{if } i=j-1,\\
z^{1/2}q^{(i-2)/2}& \text{if } i=j+1,\\
0&\text{otherwise.}
\end{cases}
\end{equation*}

Lehmer \cite{Lehmer} has computed the limit for $n\to \infty$ of the determinant of the matrix $M(n)$. Ekhad and Zeilberger \cite{EZ} have generalized this result by computing the determinant of the finite matrix $M(n)$. Furthermore, a lively account of how modern computer algebra leads to a solution   was given. Most prominently, the celebrated $q$-Zeilberger algorithm \cite{PWZ} and creative guessing were used. 

In this note, the determinant in question is obtained by computing the LU-decomposition $LU=M$. This is done with a computer, and the exact form of $L$ and $U$ is obtained by guessing. A proof that this is indeed the LU-decomposition is then a routine calculation. From it, the determinant in question is computed by multiplying the diagonal  elements of the matrix $U$. By telescoping, the final result is then quite attractive, as already stated and proved by Ekhad and Zeilberger \cite{EZ}.

We hope that this little contribution will be a welcome addition to the rekindled interest in Lehmer's tridiagonal determinant.

\medskip

We use standard notation~\cite{Andrews76}: $(x;q)_n=(1-x)(1-xq)\dots(1-xq^{n-1})$, and the Gaussian $q$-binomial coefficients
${\gauss nk}=\frac{(q;q)_{n}}{(q;q)_{k}(q;q)_{n-k}}$

\section{The LU-decomposition of $M$}

Let
\begin{equation*}
\lambda(j):=\sum_{0\le k\le j/2}\gauss{j-k}{k}(-1)^kq^{k(k-1)}z^k.
\end{equation*}
It follows from the basic recursion of the Gaussian $q$-binomial coefficients \cite{Andrews76} that
\begin{equation}\label{recu}
\lambda(j)=\lambda(j-1)-zq^{j-2}\lambda(j-2).
\end{equation}

Then we have 
\begin{equation*}
U_{j,j}=\frac{\lambda(j)}{\lambda(j-1)},\quad\text{\quad}
U_{j,j+1}=z^{1/2}q^{(j-1)/2},
\end{equation*}
and all other entries in the $U$-matrix are zero. Further,
\begin{equation*}
L_{j,j}=1,\quad\text{\quad}
L_{j+1,j}=z^{1/2}q^{(j-1)/2}\frac{\lambda(j-1)}{\lambda(j)},
\end{equation*}
and all other entries in the $L$-matrix are zero.

The typical element of the product $(LU)_{i,j}$, that is
\begin{equation*}
\sum_{1\le k\le n}L_{i,k}U_{k,j}
\end{equation*}
is almost always zero; the exceptions are as follows: If $i=j$,  then we get
\begin{equation*}
L_{j,j}U_{j,j}+L_{j,j-1}U_{j-1,j}=\frac{\lambda(j)+zq^{j-2}\lambda(j-2)}{\lambda(j-1)}=1,
\end{equation*}
because of the above recursion (\ref{recu}). If $i=j-1$,  then we get
\begin{equation*}
L_{j-1,j-1}U_{j-1,j}+L_{j-1,j-2}U_{j-2,j}=z^{1/2}q^{(j-2)/2},
\end{equation*}
and if $i=j+1$,  then we get
\begin{equation*}
L_{j+1,j+1}U_{j+1,j}+L_{j+1,j}U_{j,j}=z^{1/2}q^{(j-1)/2}\frac{\lambda(j-1)}{\lambda(j)}\frac{\lambda(j)}{\lambda(j-1)}
=z^{1/2}q^{(j-1)/2}.
\end{equation*}
This proves that indeed $LU=M$. Therefore  for the determinant of the Lehmer matrix $M$
we obtain the expression
\begin{equation*}
\prod_{j=1}^n\frac{\lambda(j)}{\lambda(j-1)}=\frac{\lambda(n)}{\lambda(0)}=
\sum_{0\le k\le n/2}\gauss{n-k}{k}(-1)^kq^{k(k-1)}z^k.
\end{equation*}
This is the result posted on August 21, 2018\footnote{Today.} by Ekhad and Zeilberger~\cite{EZ}. Of course, taking the limit $n\to\infty$,  leads
to the old result by Lehmer for the determinant of the infinite matrix:
\begin{equation*}
\lim_{n\to\infty}\det(M(n))=\sum_{k\ge0}\frac{(-1)^kq^{k(k-1)}z^k}{(q;q)_k}.
\end{equation*}

\textbf{Remarks.}

1. For $q=1$, Lehmer's determinant plays a role when enumerating   lattice paths (Dyck paths) of bounded height, or 
planar trees of bounded height, see \cite{deBrKnRi72, Knuth-Selected, HHPW}.

2. Recursions as in (\ref{recu}) have been studied in \cite{Andrews-fibo, Cigler03, PaPr03} and are linked to so-called
Schur polynomials \cite{Schur17}.

\bibliographystyle{plain}
%\bibliography{lehmer}

\end{document}